\documentclass[fontsize=11pt,headings=small,headinclude=true,footinclude=true, 
 headings=optiontohead,toc=listof,DIV=22,listof=leveldown, headsepline , letterpaper]{scrartcl} 
\usepackage[utf8]{inputenc}
\usepackage[T1]{fontenc}
\usepackage[ngerman, english]{babel}

\newcommand{\dptitlea}{An Integration--Annihilator method}
\newcommand{\dptitleb}{for analytical solutions of Partial Differential Equations}
\newcommand{\dptitle}{\dptitlea \\ \dptitleb}
\newcommand{\dptitleclean}{\dptitlea~\dptitleb}

\newcommand{\dpautoren}{%
Oliver Richters$^1$ \href{https://orcid.org/0000-0001-8253-4716}{\includegraphics[width=1em]{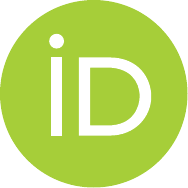}},
Erhard Glötzl$^2$ \href{https://orcid.org/0000-0002-3092-8243}{\includegraphics[width=1em]{orcid_logo.pdf}}%
}
\newcommand{\dpautorenclean}{Oliver Richters, Erhard Glötzl}

\newcommand{\dpaffiliation}{\small 1: Potsdam Institute for Climate Impact Research, Potsdam, Germany. \\ 2: Institute of Physical Chemistry, Johannes Kepler University Linz, Austria.}

\usepackage{amsmath}
\usepackage{amssymb,amsthm,amsfonts,textcomp}
\usepackage{color, transparent}
\usepackage{dpfloat}
\usepackage{array}
\usepackage{txfonts}
\usepackage{booktabs, multirow}
\usepackage{paralist}
\usepackage{hhline}
\usepackage{graphicx}
\usepackage{hanging}
\usepackage{bm} \usepackage{bbold}
\usepackage[figuresleft]{rotating}
\usepackage{acronym}
\usepackage{microtype}
\usepackage[hang]{footmisc}
\usepackage{changepage} 
\usepackage{appendix}
\usepackage{csquotes}
\usepackage{pdfpages}
\usepackage{xcolor}
\usepackage{refcount}

\makeatletter
\newcommand{\extractnumber}[1]{\expandafter\@secondoftwo\csname r@#1\endcsname}
\makeatother

\usepackage{pifont} 

\usepackage[backend=biber, natbib=true, style=authoryear-comp, sorting=nyt, url=true, isbn=false, doi=true, eprint=true, maxcitenames=2, citetracker=true, uniquename=false, labeldate=year]{biblatex}

\DeclareSourcemap{
 \maps{
  \map{
   \step[fieldsource=language, fieldset=langid, origfieldval, final]
   \step[fieldset=language, null]
  }
 }
}

\AtEveryBibitem{%
   \clearfield{day}%
   \clearfield{month}
   \clearfield{endday}%
   \clearfield{endmonth}%
}

\usepackage[inner=3.2cm,outer=3.2cm,top=1.2cm,bottom=0.8cm,includeheadfoot , heightrounded]{geometry}

\usepackage{chngcntr}

\makeatletter
\newcommand\arraybslash{\let\\\@arraycr}
\makeatother

\makeatletter
\g@addto@macro\UrlBreaks{\do\*\do\~\do\'\do\"\do\a\do\b\do\c\do\d\do\e\do\f\do\g\do\h\do\i\do\j\do\k\do%
\l\do\m\do\n\do\o\do\p\do\q\do\r\do\s\do\t\do\u\do\v\do\w\do\x\do\y\do\z\do\&\do\1\do\2\do\3\do\4\do\5\do\6\do\7\do\8\do\9\do\0\do\.}
\makeatother

\usepackage{scrlayer-scrpage}

\clearscrheadfoot
\ihead{\pagemark}
\ohead{\scriptsize Richters, Glötzl: \dptitleclean}
\let\OLDthebibliography\thebibliography
\renewcommand\thebibliography[1]{
 \OLDthebibliography{#1}
 \setlength{\parskip}{0pt}
 \setlength{\itemsep}{0pt plus 0.3ex}
}

\KOMAoption{listof}{leveldown}

\usepackage{ragged2e}
\newcolumntype{L}[1]{>{\raggedright\arraybackslash}p{#1}} 
\newcolumntype{C}[1]{>{\centering\arraybackslash}p{#1}} 
\newcolumntype{R}[1]{>{\raggedleft\arraybackslash}p{#1}} 

\DeclareMathOperator{\A}{A}
\DeclareMathOperator{\B}{B}
\DeclareMathOperator{\D}{D}
\DeclareMathOperator{\grad}{grad}
\DeclareMathOperator{\Div}{div}

\DeclareMathOperator{\ROT}{ROT}

\newcommand{\uphi}{u}
\newcommand{\x}{\bm x}

\theoremstyle{plain}
\newtheorem{theorem}{Theorem}[section]

\theoremstyle{definition}
\newtheorem{corol}[theorem]{Corollary}
\newtheorem{defin}[theorem]{Definition}
\newtheorem{lemmarev}[theorem]{Lemma}

\newcommand{\todone}[1]{}

\usepackage[hidelinks, breaklinks=true, pdfauthor={\dpautorenclean}, pdftitle={\dptitleclean}, unicode]{hyperref}

\usepackage{tikz}
\usetikzlibrary{shapes.geometric, arrows}
\tikzstyle{startstop} = [rectangle, rounded corners, minimum width=2.8cm, minimum height=1cm,text centered, draw=black, fill=red!30]
\tikzstyle{decision} = [rectangle, rounded corners, minimum width=2.8cm, minimum height=1.2cm, text centered, draw=black, fill=orange!30]
\tikzstyle{process} = [rectangle, minimum width=2.8cm, minimum height=3.2cm, text centered, draw=black, fill=green!30, execute at begin node=]
\tikzstyle{processslim} = [rectangle, minimum width=2.8cm, minimum height=1.2cm, text centered, draw=black, fill=green!30, execute at begin node=]
\tikzstyle{arrow} = [thick,->,>=stealth]

\begin{document}

\selectlanguage{english}

\thispagestyle{scrplain}

\begin{center}

 { \Large { \bfseries \sffamily \dptitle \\ \bigskip  \par \par } \vspace{1em} {\large \dpautoren \par \vspace{1em} \normalsize \dpaffiliation \par \vspace{1em} May 15, 2025
\ -- Version 1 } }
\end{center}
\vspace{1.3em}
\begin{addmargin}{0.05\textwidth}

\textbf{Abstract:}
We present a novel method to derive particular solutions for partial differential equations of the form $(\A + \B)^k Q(\x) = q(\x)$, with $\A$ and $\B$ being linear differential operators with constant coefficients, $k$ an integer, and $Q$ and $q$ sufficiently smooth functions.
The approach requires that a function $W$ and an integer $\lambda$ can be found with the following two conditions: $q$ can be integrated with respect to $\A$ such that $\A^{\lambda + k} W(\x) = q(\x)$, and $\B^{\lambda + 1}$ annihilates $W$ such that $\B^{\lambda + 1} W(\x) = 0$.

Applications include the Poisson equation $\Delta Q(\x) = q(\x)$, the inhomogeneous polyharmonic equation $\Delta^k Q(\x) = q(\x)$, the Helmholtz equation $(\Delta + \nu) Q(\x) = q(\x)$ and the wave equation $\Box Q(\x) = q(\x)$.
We show how solving the Poisson equation allows to derive the Helmholtz decomposition that splits a sufficiently smooth vector field into a gradient field and a divergence-free rotation field.

\bigskip


\bigskip

\noindent\begin{minipage}[t]{0.84\linewidth}
\textbf{Licence:} Creative-Commons \href{http://creativecommons.org/licenses/by-nc-nd/4.0/}{CC-BY-NC-ND 4.0}. \end{minipage} 
\begin{minipage}[t]{0.15\linewidth}\vspace{-\ht\strutbox}
\includegraphics[width=\columnwidth]{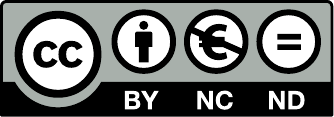}\end{minipage}

\bigskip


\end{addmargin}

\bigskip

\pagestyle{scrheadings}


\section{Introduction}

Several multivariate partial differential equations that are important in theoretical physics and applied mathematics have the functional form $\D^k Q(\x) = q(\x)$ where the differential operator $\D$ is a sum of two or more other linear differential operators with constant coefficients.
Here, $\x$ is a vector in $\mathbb{R}^n$, $k$ an integer, and $q, Q$ are sufficiently smooth scalar functions or vector components.

In this paper, we present an integration--annihilator method that allows to derive a particular solution by splitting the operator $\D$ into two components $\D = \A + \B$.
For the inhomogeneity $q(\x)$, an \emph{integration} is performed with respect to $\A$ by finding a function $W$ such that $\A^{\lambda + k} W(\x) = q(\x)$.
The integer $\lambda$ has to be chosen such that $\B^{\lambda + 1}$ is an \emph{annihilator} of $W$, meaning that $\B^{\lambda + 1} W(\x) = 0$.
A particular solution of the partial differential equation can then be expressed by a weighted sum of all possible ways to apply $\A$ and $\B$ a total of $\lambda$ times to $W$.
The weights will be chosen such that when $\D^k$ is applied, many terms cancel each other out or are omitted due to the annihilator condition.

The paper proceeds as follows:
Section~\ref{sec_annihilator} summarizes existing approaches to use annihilators for solving ordinary and partial differential equations.
Section~\ref{sec_theorem_AB} presents the main theorems which are then applied in Section~\ref{sec_appl} to the following examples of differential equations used in physics:
\begin{itemize} \itemsep0em 
\item Sec.\ \ref{sec_genpoly}: A generalized polyharmonic equation $\tilde \Delta^k Q(\x) = q(\x)$ with $\tilde \Delta = \sum_{i=1}^n \omega_i \partial_i^2$ using the example of the wave equation,
\item Sec.\ \ref{sec_poisson}: Poisson's equation $\Delta Q(\x) = q(\x)$ with $\Delta = \sum_{i=1}^n \partial_i^2$ being the Laplace operator,
\item Sec.\ \ref{sec_polyharmonic}: the inhomogeneous polyharmonic equation $\Delta^k Q(\x) = q(\x)$ (called biharmonic for $k = 2$),
\item Sec.\ \ref{sec_helmholtz}: a generalized Helmholtz equation $(\tilde \Delta^j + \nu)^k Q(\x) = q(\x)$, with the original Helmholtz equation $(\Delta + \nu) Q(\x) = q(\x)$ as a special case,
\item Sec.\ \ref{sec_hd}: the Helmholtz decomposition that splits a sufficiently smooth vector field into a gradient field and a divergence-free rotation field.
\end{itemize}
In each case, we suggest methods to choose $\A$, $\B$, $\lambda$ and $W$
to derive closed-form solutions for several functions that typically consist of polynomials, exponential functions, trigonometric functions and their products and sums.
Section~\ref{sec_conclusion} discusses our results.
The Appendix presents 3 examples and a Mathematica worksheet \citep{richters_2025_annihilatorzenodo} allows to automatically apply our theorems to further functions.

\section{Notation}
\label{sec_notation}

We denote the partial derivative as $\partial_{x_j}$, the k-th partial derivative as $\partial_{x_j}^k$, the generalized Laplace operator as $\tilde \Delta = \sum_{i=1}^n \omega_i \partial_{x_i}^2$ with $\omega_i \neq 0$ real weights, and the Laplace operator $\Delta$ as a special case for $w_i = 1$.
The generalized Laplacian to the power of $p$ is written as $\tilde \Delta^p$.
We define the incomplete generalized Laplacian (and, analogously, the incomplete Laplacian $\Delta_{\setminus k}$) and its $p$-th power as:
\begin{align}
\tilde \Delta_{\setminus k} f \coloneqq \sum_{i \neq k} \omega_i \partial_{x_i}^2 f \quad \text{ and } \quad \tilde \Delta_{\setminus k}^p f \coloneqq \tilde \Delta_{\setminus k} (\tilde \Delta_{\setminus k}^{p-1} f).
\end{align}
We denote the antiderivative of a scalar field $f$ with respect to $x_j$ as:
\begin{align}
\mathcal{A}_{x_j} f(x) \coloneqq& \int_{x_0}^{x_j} f(\xi) d\xi_j,
\end{align}
where $x_0$ can be freely chosen.
The $p$-th antiderivative of a scalar $f$ with respect to $x_j$ is given by the Cauchy formula for repeated integration or the Riemann–Liouville integral \citep{riesz_integrale_1949, cauchy_trente-cinquieme_1823}:
\begin{align}
\mathcal{A}_{x_j}^p f(x) \coloneqq&\, \mathcal{A}_{x_j}^{p-1} \mathcal{A}_{x_j} f(x) = \frac{1}{(p-1)!} \int_{x_0}^{x_j} (x_j - t_j)^{p-1} f(t) dt_j, 
\end{align}
By convention, $\partial_{x_j} \mathcal{A}_{x_j} f = \partial_{x_j}^p \mathcal{A}_{x_j}^p f = \mathcal{A}_{x_j}^0 f = \partial_{x_j}^0 f = \Delta^0 f = \Delta_{\setminus k}^0 f = \tilde \Delta^0 f = \tilde \Delta_{\setminus k}^0 f = f$.

We denote an $n$-dimensional vector field $\bm{f} = [f_k;\,_{1 \leq k \leq n}]$ with bold, italic font and use single indices for its components $f_k$.
We denote an $n \times n$ matrix-valued map $\mathbf{L} = \llbracket L_{ij} \rrbracket = \llbracket L_{ij}; \,_{1 \leq i \leq n, 1 \leq j \leq n} \rrbracket$ with bold, non-italic font.
For vector fields and matrices, the definitions for derivatives and integrals above are applied component-wise.

For $h \in \mathbb{R}$, $\lceil h \rceil$ is the least integer greater than or equal to $h$ (ceiling function).

\section{Using annihilators to solve differential equations}
\label{sec_annihilator}

The \emph{method of undetermined coefficients via annihilation}%
\footnote{The approach summarized here is often simply called \emph{annihilator method}, but because the last step uses the \emph{method of undetermined coefficients}, some textbooks use these terms as synonyms (\citealp[116]{bohner_dynamic_2012}; \citealp[234]{dharmaiah_introduction_2012}). One possible reason is that the name \emph{annihilator method} was (according to a google books and scholar search) first used in the 1960ies \citep{brauer_elementary_1968, apostol_calculus_1967}, while the concepts to annihilate the non-homogeneous term were developed much earlier, partially using the \emph{undetermined coefficients} name (for example \citealp{conkwright_method_1934}; \citealp[132--137]{jordan_cours_1915}).
The \emph{method of undetermined coefficients} based on educated guesses of the correct functional form goes back even further and was used for example by Euler to derive his lunar theories in the 1740ies \citep{verdun_leonhard_2013}.
Today, tables of typical functions simplify this approach (\citealp[88-92]{derrick_elementary_1976}; \citealp[234]{dharmaiah_introduction_2012}), and these tables can be derived and justified using the annihilator method (\citealp[264]{nagle_fundamentals_1993}; \citealp[67]{sanchez_ordinary_2002}).
To avoid a confusion of terms, we consider it convincing to call this systematic approach the \emph{method of undetermined coefficients via annihilation} as done by \citet[264--274]{wirkus_course_2014}.}
is a technique to determine a particular solution for some \emph{non-homogeneous} ordinary differential equations $\D(y(x)) = q(x)$, when $q$ is a solution of some other \emph{homogeneous} differential equation with constant coefficients $\B(q(x)) = 0$.
The latter formula is the reason why the differential operator $\B$ is said ``to annihilate'' $q(x)$ or ``to be an annihilator of'' $q(x)$ \citep[197]{farlow_introduction_2012}.
Tables exist in order to derive $\B$, for example \citep[82]{joyner_introduction_2012}:
\begin{align}
q(x) &= x^k e^{\alpha x} && \Rightarrow \B = (\partial_x - \alpha)^{k+1} \\
q(x) &= x^k e^{\alpha x} \cos(\beta x) && \Rightarrow \B = (\partial_x^2 - 2\alpha \partial_x + \alpha^2 + \beta^2)^{k+1}
\end{align}
If $\B(q(x)) = 0$, then for every solution $y(x)$ that satisfies $\D(y(x)) = q(x)$, clearly $\B(\D(y(x))) = \B(q(x)) = 0$ and therefore $y$ is a solution of the homogeneous equation $\B(\D(y(x))) = 0$.
Thus, if the homogeneous equation $\B(\D(y(x))) = 0$ can be solved, for example using the technique of the \emph{characteristic equation}, its set of solutions with yet undetermined coefficients can be used as an ansatz for $\D(y(x)) = q(x)$.
Then,  by comparing $\D(y(x))$ with $q(x)$, the coefficients that lead to the final non-homogeneous solution can be derived (\citealp[90--93]{coddington_introduction_2012}; \citealp[264--274]{wirkus_course_2014}).%
\footnote{%
  As an example, let us look at the differential equation $\D(y(x)) = \partial_x^2 y - 3 \partial_x y + 2y = x^2$.
  Since $x^2$ is a solution of $\B(y(x)) = \partial_x^3 y = 0$, every solution of $\D(y(x)) = x^2$ is a solution of
  $\B(\D(y(x))) = \partial_x^5 y - 3 \partial_x^4 y + 2 \partial_x^3 y = 0$. The characteristic equation of this equation is
  $r^3(r^2 - 3 r +2)$, the product of the characteristic polynomials for $\D$ and $\B$. The roots are 0, 0, 0, 1, 2, and hence
  $y(x)$ must have the form $y(x) = c_0 + c_1 x + c_2 x^2 + c_3 e^x + c_4 e^{2x}$. Obviously,
  $c_3 e^x + c_4 e^{2x}$ solves $\D(y(x)) = 0$, and comparison of the other coefficients leads to $c_2 = 1/2$,
  $c_1 = 3/2$ and $c_0 = 7/4$ which is a particular solution of $\D(y(x)) = x^2$
  \citep[90--93]{coddington_introduction_2012}.%
}

Using annihilators for ordinary differential equations is a systematic and efficient method when applicable, but is limited to $q(x)$ having the right functional form such that an annihilator exists \citep[166]{apostol_calculus_1967}.

Adapting the method described above to \emph{partial} differential equations is straightforward.
For example, \citet{golberg_annihilator_1999} uses it to approximate a function using basis functions that can be integrated analytically, and then a comparison of coefficients is performed.
The disadvantage of this approach is that for each different inhomogeneity, a different annihilator is used, and the full process has to be repeated.

In the following, we describe a new approach that differs in two aspects:
First, it is restricted to finding the particular solution, separating the homogeneous and the inhomogeneous problem.
Second, it automatically generates the correct coefficients, avoiding a comparison of coefficients.

\section{An Integration--Annihilator method}

\label{sec_theorem_AB}
In this section, we present a novel approach to derive a closed-form particular solution for $(\A + \B)^n Q = q$ with $\A$ and $\B$ linear differential operators with constant coefficients.
It combines integration with respect to $\A$ to derive a function $W$ with the repeated application of $\B$ to annihilate $W$. The restriction to \emph{linear with constant coefficients} guarantees that $\A$, $\B$ and $\D$ commute.

Strategies on how to determine $\A$, $\B$, $\lambda$ and $W$ for various partial differential equations will be presented in Section~\ref{sec_appl}.
Concrete examples of the use of the following theorems can be found in the Appendix.

\begin{theorem} \label{theorem_AB} 
Let $q \in C^1(\mathbb{R}^n, \mathbb{R})$ be a scalar function, let $\A$, $\B$ and $\D$ be some linear differential operators with constant coefficients and $\D = \A + \B$, and let there exist two natural
numbers $\lambda$ and $k$ and a sufficiently smooth function $W$ such that:
\begin{align}
\A^{\lambda + k} W(\x) &= q(\x) && \text{(integration condition)}, \label{eq_bedingung_ABop1} \\
\B^{\lambda + 1} W(\x) &= 0 && \text{(annihilator condition)}. \label{eq_bedingung_ABop2}
\end{align}
Then, $Q$ as defined below satisfies $\D^k Q = q$:
\begin{align}
Q(\x) &\coloneqq \sum_{p=0}^{\lambda}  (-1)^{p} \tbinom{k+p-1}{p} \A^{\lambda - p} \B^p W(\x). \label{eq_Q_AB}
\end{align}
\end{theorem}

\begin{proof}
The proof follows by complete induction. We will show that for $0 \leq j \leq k$, $\D^j Q(\x)$ is given by:
\begin{align}
 \D^j Q(\x) &= \A^{\lambda + j} W(\x) + \sum_{p=1}^{\lambda} (-1)^p \tbinom{k+p-j-1}{p} \A^{\lambda-p+j} \B^p W(\x). \label{eq_Q_AB_proof1}
\end{align}

Base case for $j = 0$: the fact that $\D^0 Q(\x) = Q(\x)$ follows by identifying the first term in Eq.~\eqref{eq_Q_AB_proof1} with the $p=0$ term from the sum in Eq.~\eqref{eq_Q_AB} because $(-1)^0 \tbinom{k-1}{0} = 1$, while the terms for $p > 0$ are identical.

Induction step: If Eq.~\eqref{eq_Q_AB_proof1} is true for $j < k$ (induction hypothesis), we will show that it is true also for $j+1$.
First, as $(-1)^0 \binom{k+p-j-1}{0} = 1$, we include the first term of Eq.~\eqref{eq_Q_AB_proof1} into the sum for $p=0$:
 \begin{align}
  \D^j Q(\x) &= \sum_{p=0}^{\lambda} (-1)^p \tbinom{k+p-j-1}{p} \A^{\lambda-p+j} \B^p W(\x). \\
  \intertext{Now, we apply another operator $\D$ and use that $\A$ and $\B$ commute:}
  \D^{j+1} Q(\x)
   &= \D \D^j Q(\x) = (\A + \B) \D^j Q(\x) \\
   &= \sum_{p=0}^{\lambda}   (-1)^p \tbinom{k+p-j-1}{p} \A^{\lambda-p+j+1} \B^p     W(\x) 
                   + \sum_{p=0}^{\lambda}   (-1)^p \tbinom{k+p-j-1}{p}   \A^{\lambda-p+j}  \B^{p+1} W(\x), \\
   \intertext{now shift the $p$ index by one in the second sum}
   &= \sum_{p=0}^{\lambda}   (-1)^p \tbinom{k+p-j-1}{p}   \A^{\lambda-p+j+1} \B^p W(\x) 
                  - \sum_{p=1}^{\lambda+1} (-1)^p \tbinom{k+p-j-2}{p-1} \A^{\lambda-p+j+1} \B^p W(\x), \\
   \intertext{now the $p=0$ term of the first sum can be taken out of the sum, and the $\lambda +1$ term of the second sum is zero because $\B^{\lambda +1} W(\x) = 0$ by definition. At the same time, as it is well-known that $\tbinom{a}{b} -  \tbinom{a-1}{b-1} = \tbinom{a-1}{b}$ for all $a, b$, all terms with equal $p$ can be summed together to yield:}
   &= \A^{\lambda + j +1} W(\x)
      + \sum_{p=1}^{\lambda} (-1)^p \tbinom{k+p-(j+1)-1}{p} \A^{\lambda-p+j+1} \B^p W(\x)
  \end{align}
which is exactly Eq.~\eqref{eq_Q_AB_proof1} for $j+1$.

Conclusion: Since both the base case and the induction step have been proven as true, by mathematical induction the statement in Eq.~\eqref{eq_Q_AB_proof1} holds for every integer $j \leq k$.

Now, setting in Eq.~\eqref{eq_Q_AB_proof1} $j = k$, then $k+p-j-1$ becomes $p-1$, and $\tbinom{p-1}{p} = 0$ for all $p$, thus the sum can be dropped, so only the first term is non-zero and yields, using Eq.~\eqref{eq_bedingung_ABop1}:
\begin{align}
 \D^k Q(\x) &= \A^{\lambda + k} W(\x) = q(\x).
\end{align}
This finalizes the proof of Theorem~\ref{theorem_AB}.
\end{proof}

In the special case $k=1$ and therefore $\D Q(\x) = q(\x)$, we will now relax the condition that $\B^{\lambda + 1}$ annihilates $q(\x)$ to the condition that $\B^{\lambda+1}$ yields some multiple of $q(\x)$.

\begin{theorem} \label{theorem_AB1}

Let $q \in C^1(\mathbb{R}^n, \mathbb{R})$ be a scalar function, let $\A$, $\B$ and $\D$ be some linear differential operators with constant coefficients and $\D = \A + \B$, and let there exist a non-zero real number $\uphi$, a natural
number $\lambda$ and a sufficiently smooth function $W$ such that:
\begin{align}
\A^{\lambda + 1} W(\x) &= q(\x) && \text{(integration condition)}, \label{eq_bedingung_AB1} \\
\B^{\lambda + 1} W(\x) &= (-1)^{\lambda} (\uphi - 1) q(\x) && \text{(relaxed annihilator condition)}. \label{eq_bedingung_AB2}
\end{align}
Then, $Q$ as defined below satisfies $\D Q = q$:
\begin{align}
Q(\x) &\coloneqq \sum_{p=0}^{\lambda} \frac{(-1)^{p}}{\uphi} \A^{\lambda - p} \B^p W(\x). \label{eq_Q_AB1}
\end{align}
\end{theorem}

\begin{proof}
We apply another operator $\D$
\begin{align}
\D Q(\x) &= \D \sum_{p=0}^\lambda \frac{(-1)^{p}}{\uphi} \A^{\lambda - p} \B^p W(\x), \\
\intertext{replace $\D$ with $\A + \B$ and use that $\A$ and $\B$ commute}
&= \frac{1}{\uphi} \sum_{p=0}^{\lambda} \left[ (-1)^p \A^{\lambda - p + 1} \B^p W(\x) + (-1)^p \A^{\lambda - p} \B^{p+1} W(\x) \right], \\
\intertext{shift the sum index $p$ by 1 in the second sum}
&= \frac{1}{\uphi} \left[ \sum_{p=0}^{\lambda} (-1)^p \A^{\lambda-p+1} \B^p W(\x) - \sum_{p=1}^{\lambda+1} (-1)^p \A^{\lambda-p+1} \B^p W(\x) \right], \\
\intertext{and cancel the terms with equal index $p$, keeping only the terms with $p=0$ and $p=\lambda+1$}
&= \frac{1}{\uphi} \left[ \A^{\lambda+1} W(\x) + (-1)^{\lambda} \B^{\lambda + 1} W(\x) \right] = \frac{1}{\uphi} \left[ q(\x) + (\uphi - 1) q(\x) \right]. \\
&= q(\x).
\end{align}
This finalizes the proof of Theorem~\ref{theorem_AB1}.
\end{proof}

The following trivial statement will be used to cover some special cases later:

\begin{lemmarev} \label{lemma_sin}
Let $\D$ be a differential operator and let $q \in C^\infty(\mathbb{R}^n, \mathbb{R})$ be a scalar function.
If $\D q(\x) = v \cdot q(\x)$ with $v$ being a non-zero real constant, then $Q(\x) = \frac{1}{v^k} q(\x)$ is a solution of $\D^k Q(\x) = q(x)$.
\end{lemmarev}
\begin{proof}
The proof follows immediately from applying $\D$ $k$-times to $Q$.
\end{proof}


\section{Applications}
\label{sec_appl}

In this section, the theorems from the previous section are applied to some important partial differential equations from physics:
First, we show a rather general proof for a generalized Laplace operator $\tilde \Delta^k Q(\x) = q(\x)$ using the example of the inhomogeneous wave equation.
Then, we show special cases, namely Poisson's equation $\Delta Q(\x) = q(\x)$ (Section~\ref{sec_poisson}), and the polyharmonic equation $\Delta^k Q(\x) = q(\x)$
(Section~\ref{sec_polyharmonic}) using the standard Laplace operator.
The flowchart in Figure~\ref{fig_flowchart} shows how to select the proper theorem.
Additional, we provide solutions for a generalized Helmholtz equation (Section~\ref{sec_helmholtz}) and the Helmholtz decomposition (Section~\ref{sec_hd}).
A \emph{Mathematica} worksheet implementing the solution strategies can be found in \citet{richters_2025_annihilatorzenodo} and at \url{https://oliver-richters.de/annihilator}.

\subsection{Generalized Laplace operator / inhomogeneous wave equation}
\label{sec_genpoly}

\begin{corol} 
\label{theorem_WPH}

Let $\tilde \Delta = \sum_{i=1}^n \omega_i \partial_{x_i}^2$ be a generalized Laplace operator with $\omega_i \neq 0$ and let $q \in C^1(\mathbb{R}^n, \mathbb{R})$ be a scalar function such that there $\exists m \leq n$ with $q(\x) = f(x_m) \cdot \x_{i \neq m}^{\beta}$ with multi-index $\beta \in \mathbb{N}_0^{n-1}$, thus the inhomogeneity can be separated into a monomial and a function that depends on a single coordinate $x_m$. Then, choose $\lambda = \lceil \tfrac{|\beta| - 1}{2} \rceil$ and $W(\x) = \omega_m^{-\lambda - k} \mathcal{A}_{x_m}^{2 \lambda + 2 k} q(\x) = \omega_m^{-\lambda - k} \x_{i \neq m}^{\beta} \mathcal{A}_{x_m}^{2 \lambda + 2 k} f(x_m)$.
The solution $Q$ as defined below satisfies $\tilde \Delta^k Q = q$:
\begin{align}
Q(\x) &\coloneqq \sum_{p=0}^{\lambda} (-1)^{p} \tbinom{k+p-1}{p} \omega_m^{\lambda - p} \partial_{x_m}^{2\lambda - 2p} \tilde \Delta_{\setminus m}^{p} W(\x) = \sum_{p=0}^{\lambda} (-1)^p \tbinom{k+p-1}{p} \omega_m^{- p - k} \mathcal{A}_{x_m}^{2 k + 2 p} \tilde \Delta_{\setminus m}^{p} q(\x). 
\end{align}
For $k=1$, this simplifies to 
\begin{align}
Q(\x) &\coloneqq \sum_{p=0}^{\lambda} (-1)^{p} \omega_m^{\lambda - p} \partial_{x_m}^{2\lambda - 2p} \tilde \Delta_{\setminus m}^{p} W(\x) = \sum_{p=0}^{\lambda} (-1)^p \omega_m^{- p - 1} \mathcal{A}_{x_m}^{2 + 2 p} \tilde \Delta_{\setminus m}^{p} q(\x).
\end{align}
\end{corol}
\begin{proof}
The proof follows from Theorem~\ref{theorem_AB} by setting $\A = \omega_m \partial_{x_m}^2$, $\B = \tilde \Delta_{\setminus m}$ and $\D = \tilde \Delta$.
The annihilator condition $\tilde \Delta_{\setminus m}^{\lambda + 1} W(\x) = 0$ is satisfied as at least $|\beta| + 1$ derivatives are applied to $x_{i \neq m}^\beta$, yielding zero.
The integration condition is satisfied as $\A^{\lambda + k} W(\x) = \omega_m^{\lambda + k} \partial_{x_m}^{2\lambda + 2k} W(\x) = q(\x)$.
\end{proof}

Keep in mind that the definition of $\mathcal{A}_{x_m}$ allows to choose the lower boundary of the integral in such a way that the calculation is as simple as possible.

An application of Corollary~\ref{theorem_WPH} with $k = 1$ is the inhomogeneous wave equation that describes the propagation of waves in the presence of external sources or forces with applications in acoustics, electromagnetics, and seismology.
Here, we use $t$ for the temporal dimension of the four-dimensional spacetime, and $\Box$ is D'Alembert's operator.
\begin{align}
\Box Q(\x) = \left( \partial_t^2 / c^2 - \sum_{i=1}^3 \partial_{x_i}^2 \right) Q(\x) = q(\x),
\end{align}
with $c$ the wave speed. \ref{app_wave} in the Appendix shows that for $q(\x) = t \sin(t) x_1^2 x_2$, a special solution of $\Box Q = q$ is given by $Q(\x) = c^2 x_1^2 x_2 \big(-2 \cos(t) - t \sin(t) \big)  + 2 c^4 x_2 \big(4 \cos(t) + t \sin(t) \big)$.

\subsection{Poisson's equation}
\label{sec_poisson}

Poisson's equation is an elliptic partial differential equation named after Denis \citet{poisson_memoire_1823}.
For scalar functions $q$ and $Q$, it is given by $\Delta Q(\x) = q(\x)$ with $\Delta$ the Laplace operator.
The homogeneous variant $\Delta Q(\x) =  0$ is called Laplace's equation \citep{hackbusch_poisson_2017}.

Poisson's equation has various applications in theoretical physics:
In Newtonian physics, the relation between the gravitational potential $\varphi_g$ and the mass density distribution $\rho_g$ is given by
$\Delta \varphi_g = 4\pi G \rho_g$ with $G$ the gravitational constant.
In electrostatics, the electrostatic potential $\varphi_e$ and the density distribution of electric charges $\rho_e$ is connected by $\Delta \varphi_e = -\frac{\rho_f}{\varepsilon}$, where $\varepsilon$ is the permittivity of the medium \citep{hackbusch_poisson_2017}.
For scalar functions that decay sufficiently fast at infinity, it can be solved with Kernel integrals (Green's function),
but for unboundedly growing vector fields, these integrals diverge \citep{chorin_mathematical_1990, schwarz_hodge_1995}.

\begin{corol}
\label{theorem_P} 

Let $q \in C^1(\mathbb{R}^n, \mathbb{R})$ be a scalar function and let there $\exists m$ with $q(\x) = f(x_m) \cdot \x_{i \neq m}^{\beta}$ with multi-index $\beta \in \mathbb{N}_0^{n-1}$ and $f$ any function of $x_m$ that can be integrated.
Then, $Q$ as defined below with $\lambda = \lceil \tfrac{|\beta| - 1}{2} \rceil$ and $W(\x) = \mathcal{A}_{x_m}^{2 \lambda + 2} q(\x) = \x_{i \neq m}^{\beta} \mathcal{A}_{x_m}^{2 \lambda + 2} f(x_m)$ satisfies $\Delta Q = q$:
\begin{align}
Q(\x) &\coloneqq \sum_{p=0}^{\lambda} (-1)^{p} \partial_{x_m}^{2\lambda - 2p} \Delta_{\setminus m}^{p} W(\x) = \sum_{p=0}^{\lambda} (-1)^p \mathcal{A}_{x_m}^{2 + 2 p} \Delta_{\setminus m}^{p} q(\x). \label{eq_potential_F}
\end{align}
\end{corol}
\begin{proof}
The proof follows from Corollary~\ref{theorem_WPH} by setting $\omega_i = 1$ and $k = 1$.
\end{proof}

\ref{app_P_monomial} in the Appendix illustrates how $Q(\x) = \frac{91 x_1^2 x_2^{12} - x_2^{14}}{12012}$ can be found as a particular solution of $\Delta Q(\x) = q(\x) = x_1^2 x_2^{10}$, and why choosing $m$ matters for the brevity of the result.

\subsection{Polyharmonic equation}
\label{sec_polyharmonic}

The polyharmonic equation $\Delta^k Q(\x) = q(\x)$ is used for the analysis of the statics and vibrations of flexural plates \citep{manolis_non-homogeneous_2003, providakis_dynamic_1999} and for image interpolation, inpainting and surface reconstruction \citep{saito_polyharmonic_2006, kobayashi_image_2011}.

\begin{corol} 
\label{theorem_PH}
Let $q \in C^1(\mathbb{R}^n, \mathbb{R})$ be a scalar function, such that there $\exists m$ with $q(\x) = f(x_m) \cdot \x_{i \neq m}^{\beta}$ with multi-index $\beta \in \mathbb{N}_0^{n-1}$ and $f$ any function of $x_m$ that can be integrated.
Then, $Q$ as defined below with $\lambda = \lceil \tfrac{|\beta| - 1}{2} \rceil$ and $W(\x) = \mathcal{A}_{x_m}^{2 \lambda + 2 k} q(\x) = \x_{i \neq m}^{\beta} \mathcal{A}_{x_m}^{2 \lambda + 2 k} f(x_m)$ satisfies $\Delta^k Q = q$:
\begin{align}
Q(\x) &\coloneqq \sum_{p=0}^{\lambda} (-1)^{p} \tbinom{k+p-1}{p} \partial_{x_m}^{2\lambda - 2p} \Delta_{\setminus m}^{p} W(\x) = \sum_{p=0}^{\lambda} (-1)^p \tbinom{k+p-1}{p} \mathcal{A}_{x_m}^{2 k + 2 p} \Delta_{\setminus m}^{p} q(\x). \label{eq_Qpolyharmonic}
\end{align}
\end{corol}

\begin{proof}
The proof follows from Corollary~\ref{theorem_WPH} by setting $\omega_i = 1$.
\end{proof}

\ref{app_PH_sinus} in the Appendix illustrates how $Q(\x) = (-1)^k \sin(x_2) (x_1^2 + 2k)$ can be found as a particular solution of $\Delta^k Q(\x) = q(\x) = x_1^2 \sin(x_2)$.

Note that \citet{bondarenko_operator_1984} already found a special case of Corollary~\ref{theorem_PH} for monomials and, by linearity, polynomials.
If $\Delta^m u = x^a/a!$ (in multi-index notation), then for $2q > a_1 + \ldots + a_{n-1}$, he presented this solution \citep[cited after][]{karachik_solution_2010}:
\begin{align}
 u(x) = \sum_{i = 0}^{q-1} (-1)^i \tbinom{i + m - 1}{m-1} \frac{x_n^{2i + 2 m + a_n}}{(2i + 2 m + a_n)!} \Delta_{\setminus n} x_{\setminus n}^a/a!
\end{align}

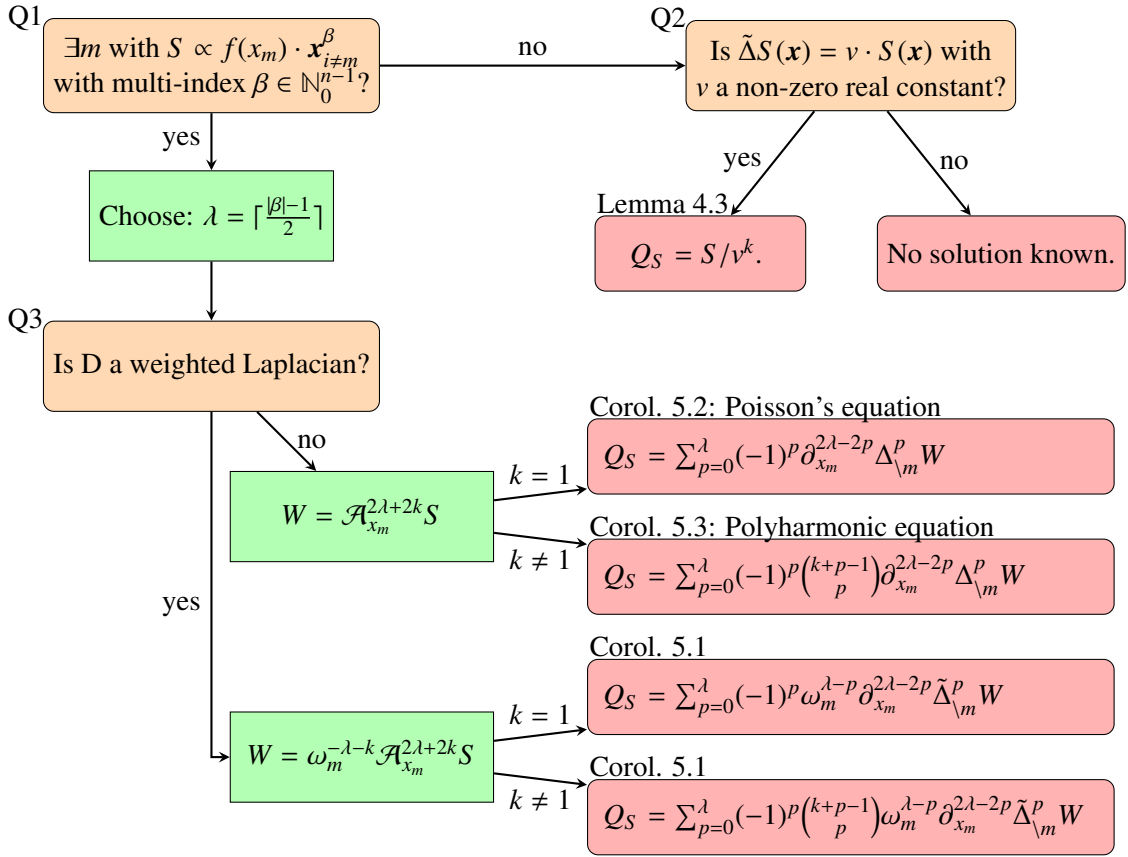
\begin{figure}[tb!]
\centering

\begin{tikzpicture}[node distance=2cm]

\newcommand{\ofx}{} 


\node[align = center] (dec_ismonomxi) [decision] {$\exists m$ with $S \propto f(x_m) \cdot \x_{i \neq m}^{\beta}$ \\ with multi-index $\beta \in \mathbb{N}_0^{n-1}$?};

\node[align = center] (dec_cossinexp) [decision, right of=dec_ismonomxi, xshift = 6.5cm] {Is $\tilde \Delta S(\x) = v \cdot S(\x)$ with \\ $v$ a non-zero real constant?};

\node[align = center] (proc_integratek) [processslim, below of = dec_ismonomxi] {Choose: $\lambda = \lceil \tfrac{|\beta| - 1}{2} \rceil$};

\node[align = center] (dec_operator) [decision, below of=proc_integratek] {Is $\D$ a weighted Laplacian?};

\draw [arrow] (dec_ismonomxi) -- node[anchor=south] {no} (dec_cossinexp);
\draw [arrow] (dec_ismonomxi) -- node[anchor=east] {yes} (proc_integratek);

\node[align = center] (proc_integrateW) [processslim, below of = dec_operator, right of = dec_operator, minimum width = 3.5cm] {$W\ofx = \mathcal{A}_{x_m}^{2 \lambda + 2 k} S\ofx$}; 

\node[align = center] (proc_integrateWw) [processslim, below of = proc_integrateW, right of = dec_operator, yshift = -3.2cm, minimum width = 3.5cm] {$W\ofx = \omega_m^{-\lambda - k} \mathcal{A}_{x_m}^{2 \lambda + 2 k} S\ofx$}; 

\node[align = center] (solutionW) [startstop, right of = proc_integrateW, yshift=0.8cm, xshift = 4.5cm, minimum width = 7cm] { %
\parbox{6.5cm}{$Q_S\ofx = \sum_{p=0}^{\lambda} (-1)^{p} \partial_{x_m}^{2\lambda - 2p} \Delta_{\setminus m}^{p} W\ofx$}
 };

\node[align = center] (solutionWk) [startstop, right of = proc_integrateW, yshift=-0.8cm, xshift = 4.5cm, minimum width = 7cm] { %
\parbox{6.5cm}{$Q_S\ofx = \sum_{p=0}^{\lambda} (-1)^{p} \tbinom{k+p-1}{p} \partial_{x_m}^{2\lambda - 2p} \Delta_{\setminus m}^{p} W\ofx$}
 };
\draw [arrow] (proc_integrateW) -- node[anchor=south] {$k=1$} (solutionW);
\draw [arrow] (proc_integrateW) -- node[anchor=north] {$k\neq 1$} (solutionWk);
 
\node[align = center] (solutionWw) [startstop, right of = proc_integrateWw, yshift=0.8cm, xshift = 4.5cm, minimum width = 7cm] { %
\parbox{6.5cm}{$Q_S\ofx = \sum_{p=0}^{\lambda} (-1)^{p} \omega_m^{\lambda - p}  \partial_{x_m}^{2\lambda - 2p} \tilde \Delta_{\setminus m}^{p} W\ofx$}
 };

\node[align = center] (solutionWwk) [startstop, right of = proc_integrateWw, yshift=-0.8cm, xshift = 4.5cm, minimum width = 7cm] { %
\parbox{6.5cm}{$Q_S\ofx = \sum_{p=0}^{\lambda} (-1)^{p} \tbinom{k+p-1}{p}  \omega_m^{\lambda - p}  \partial_{x_m}^{2\lambda - 2p} \tilde \Delta_{\setminus m}^{p} W\ofx$}
 };
\draw [arrow] (proc_integrateWw) -- node[anchor=south] {$k=1$} (solutionWw);
\draw [arrow] (proc_integrateWw) -- node[anchor=north] {$k\neq 1$} (solutionWwk);

\node[align = center] (solution2) [startstop, below of = proc_integratek, left of = dec_cossinexp, yshift=-0.5cm, xshift = 0cm] { $Q_S\ofx = S\ofx /v^k$. };

\node[align = left] (proc_noidea) [startstop, below of = proc_integratek, right of = dec_cossinexp, yshift=-0.5cm, xshift=0cm] {No solution known.};

\draw [arrow] (dec_cossinexp) -- node[anchor=west] {no} (proc_noidea);
\draw [arrow] (proc_integratek) -- (dec_operator);
\draw [arrow] (dec_cossinexp) -- node[anchor=east] {yes} (solution2);

\draw [arrow] (dec_operator) -- node[anchor=west] {no} (proc_integrateW);
\draw [arrow] (dec_operator) |- node[anchor=east, yshift=+2cm] {yes} (proc_integrateWw);


\foreach \anc/\n in {north west/Q1}{\node[anchor=\anc, xshift=-0.6cm, yshift=+0.3cm] at (dec_ismonomxi.\anc) {\n};}
\foreach \anc/\n in {north west/Q2}{\node[anchor=\anc, xshift=-0.6cm, yshift=+0.3cm] at (dec_cossinexp.\anc) {\n};}
\foreach \anc/\n in {north west/Q3}{\node[anchor=\anc, xshift=-0.6cm, yshift=+0.3cm] at (dec_operator.\anc) {\n};}

\foreach \anc/\n in {north west/Lemma \ref{lemma_sin}}{\node[anchor=\anc, xshift=-0.1cm, yshift=+0.43cm] at (solution2.\anc) {\n};}
\foreach \anc/\n in {north west/Corol. \ref{theorem_P}: Poisson's equation}{\node[anchor=\anc, xshift=-0.1cm, yshift=+0.43cm] at (solutionW.\anc) {\n};}
\foreach \anc/\n in {north west/Corol. \ref{theorem_PH}: Polyharmonic equation}{\node[anchor=\anc, xshift=-0.1cm, yshift=+0.43cm] at (solutionWk.\anc) {\n};}
\foreach \anc/\n in {north west/Corol. \ref{theorem_WPH}}{\node[anchor=\anc, xshift=-0.1cm, yshift=+0.43cm] at (solutionWwk.\anc) {\n};}
\foreach \anc/\n in {north west/Corol. \ref{theorem_WPH}}{\node[anchor=\anc, xshift=-0.1cm, yshift=+0.43cm] at (solutionWw.\anc) {\n};}

\end{tikzpicture}
\caption{\label{fig_flowchart}We suggest the following strategy to derive a particular solution for $\tilde \Delta^k Q(\x) = q(\x)$, the polyharmonic equation with a weighted, generalized Laplace operator $\D = \tilde \Delta = \sum_{i=1}^n \omega_i \partial_{x_i}^2$ (the unweighted special case for $\D = \Delta$ is the polyharmonic equation, and for $k=1$ is Poisson's equation): (1) expand $q$ into a sum of expressions $S(\x)$, (2) use the flowchart to find the appropriate solution for each expression $S(\x)$ separately, (3) sum all the resulting functions $Q_S(\x)$. \newline
In the diagram, $f(x_m)$ is any function of its argument. In Q1, if $S$ is a monomial and all $m$ satisfy the condition, we recommend to choose $m$ such that $x_m$ has the highest exponent to make $\lambda$ as small as possible and to have the fewest terms in the sum.}
\end{figure}

\subsection{Generalized Helmholtz equation}
\label{sec_helmholtz}

The Helmholtz equation is a partial differential equation first used by \citet{helmholtz_theorie_1860} to study air oscillations in tubes.
Its inhomogeneous variant has the form
\begin{align}
(\Delta + \nu) Q(\x) = q(\x),
\end{align}
representing the time-independent form of the wave equation.
As $\sin(\sqrt{\nu}x)$ solves the one-dimensional homogeneous differential equation, $\sqrt{\nu}$ is usual called the wave number.
It is crucial for studying stationary oscillating processes with applications in physics, including acoustics, electromagnetics, and quantum mechanics for problems related to wave propagation and vibration.
In the following, we study a generalized form of the Helmholtz equation, of which the original one is a special case with $j = k = 1$:
\begin{align}
(\Delta^j + \nu)^k Q(\x) = q(\x).
\end{align}
A closed-form particular solution can be obtained using the integration--annihilator method is if $q$ is a polynomial.

\begin{corol}
\label{theorem_Helmholtz}
Let $q \in C^1(\mathbb{R}^n, \mathbb{R})$ be a polynomial with highest total degree $|\beta|$.
Then, with $\lambda = \lceil \frac{|\beta|-1}{2 j} \rceil$, $Q$ as defined below satisfies $(\Delta^j + \nu)^k Q(\x) = q(\x)$:
\begin{align}
Q(\x) = \sum_{p=0}^{\lambda} (-1)^p \tbinom{k+p-1}{p} \nu^{- p - k} \Delta^{jp} q(\x). \label{eq:helmholtz}
\end{align}
\end{corol}
\begin{proof}
The proof follows from Theorem~\ref{theorem_AB}, if one chooses $\D = \Delta^j + \nu$, $\A = \nu$, $\B = \Delta^j$ and $W(\x) = \frac{1}{\nu^{\lambda + k}} q(\x)$.
Then, Eq.~\eqref{eq_Q_AB} can be simplified as follows:
\begin{align}
Q(\x) = \sum_{p=0}^{\lambda} (-1)^p \tbinom{k+p-1}{p} \nu^{\lambda - p} \Delta^{jp} W(\x) = \sum_{p=0}^{\lambda} (-1)^p \tbinom{k+p-1}{p} \nu^{- p - k} \Delta^{jp} q(\x).
\end{align}
This finalizes the proof of Corollary~\ref{theorem_Helmholtz}.
\end{proof}

For the original Helmholtz equation, Eq.~\eqref{eq:helmholtz} simplifies with $j = k = 1$ to:
\begin{align}
Q(\x) = \sum_{p=0}^{\lambda} (-1)^p \nu^{- p - 1} \Delta^{p} q(\x).
\end{align}
If, for example, $q(\x) = x_1^4 x_2^3$, then with $\lambda = 3$ and $W(\x) = \tfrac{1}{\nu^4} x_1^4 x_2^3$, a particular solution is given by:
\begin{align}
Q(\x) &= \frac{x_1^4 x_2^3}{\nu}  - \frac{12 x_1^2 x_2^3 + 6 x_1^4 x_2^1}{ \nu^2 } + \frac{24 x_2^3 + 144 x_1^2 x_2}{\nu^3} - \frac{432 x_2}{\nu^4}.
\end{align}

\subsection{Helmholtz Decomposition}
\label{sec_hd}

\begin{figure}[ht]
\hfil\includegraphics[width=0.66\columnwidth]{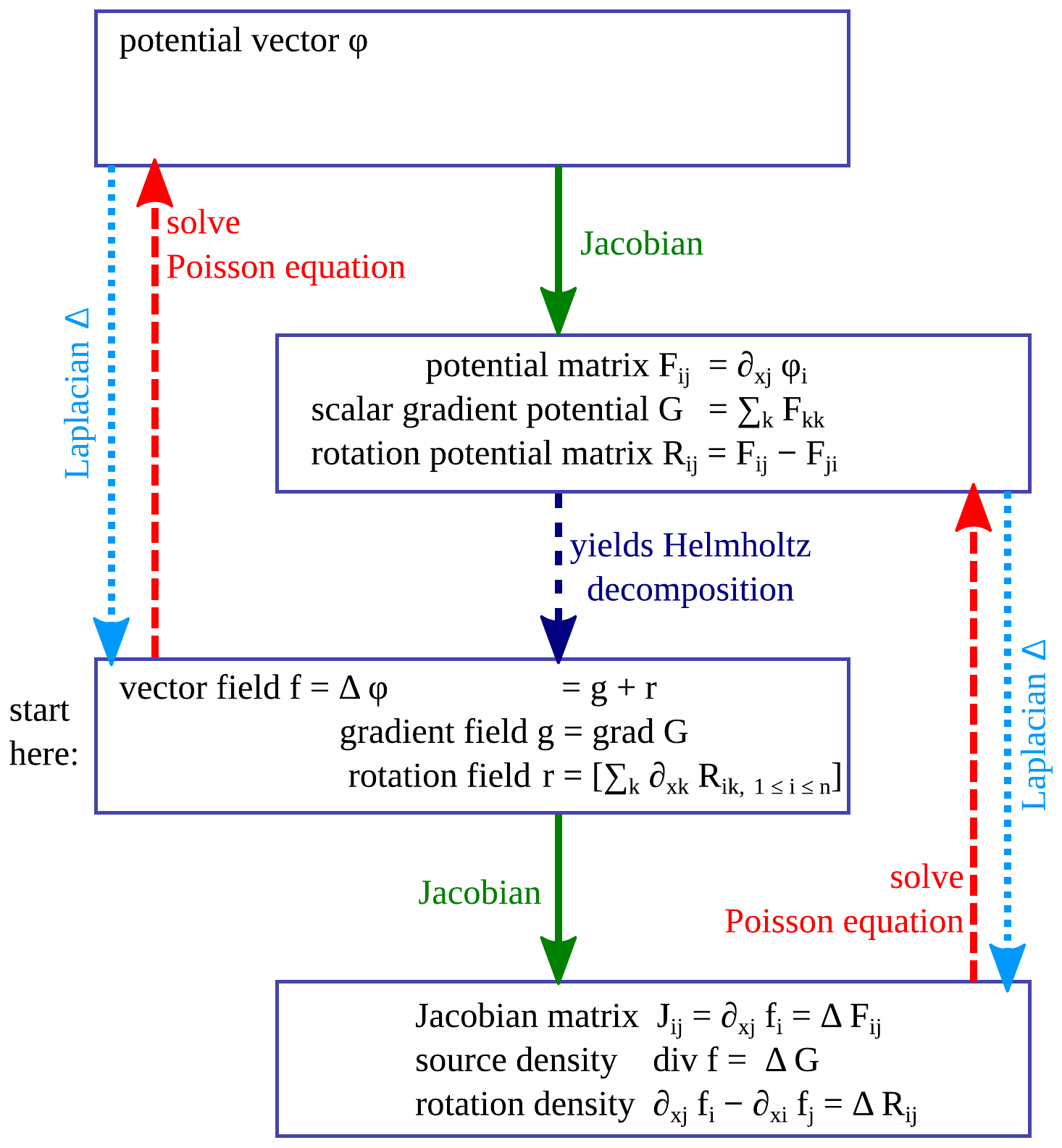}\hfil
\caption{\label{fig_zusammenhang-groessen} This figure depicts how to obtain a Helmholtz decomposition for a vector field $\bm f$ by calculating the gradient field and the rotation field from the potential matrix $\mathbf{F} = \llbracket F_{ij} \rrbracket$ (dark blue short-dashed arrow).
The variables and operators are defined in Definitions~\ref{def_helmholtz}--\ref{def_potentialmatrix}. 
This paper derives methods to solve Poisson's equation (red long-dashed arrows) which allows to determine the potential vector $\phi$ from which all other quantities can be derived.
In the traditional approach, often only in 3-dimensional space, Poisson's equation on unbounded domains only for sufficiently fast decaying fields for the source and rotation density, yielding a vector and scalar potential using the convolution of the Jacobian $\mathbf{J}$ with a solution $K$ of Laplace's equation is numerically solved.}
\end{figure}

The Helmholtz decomposition \citep{stokes_dynamical_1849, von_helmholtz_uber_1858} splits a sufficiently smooth vector field $\bm f$ into an irrotational (curl-free) `gradient field' $\bm g$ and a solenoidal (divergence-free) `rotation field' $\bm r$. It is indispensable for mathematical physics \citep{tran-cong_helmholtzs_1993, kustepeli_helmholtz_2016, dassios_uniqueness_2002, sprossig_helmholtz_2009}, but also used in animation, computer vision or robotics \citep{bhatia_helmholtz-hodge_2013}.

\begin{defin}
\label{def_helmholtz}
For a vector field $\bm{f} \in C^1(\mathbb{R}^n, \mathbb{R}^n)$, a \emph{Helmholtz decomposition} is a pair of vector fields $\bm{g} \in C^1(\mathbb{R}^n, \mathbb{R}^n)$ and $\bm{r} \in C^1(\mathbb{R}^n, \mathbb{R}^n)$ such that
$\bm{f} = \bm{g} + \bm{r}$ where $\bm{g} = \grad G$ for some function $G \in C^2(\mathbb{R}^{n}, \mathbb{R})$, and $\Div \bm{r} = 0$. The vector field $\bm{g}$ is called a \emph{gradient field} and $\bm{r}$ is called a \emph{rotation field}.
\end{defin}

If $\bm{r} = 0$, then $\bm{f}$ is a called a `curl-free' vector field and a line integral yields the gradient potential $G$.
For other cases, the usual approach is to solve Poisson's equation using Kernel integrals:\footnote{Note that in physics, $G$ is often defined with opposite sign \citep[p.~105]{duschek_grundzuge_1961}.}
\begin{align}
\Delta G(\x) &= \Div \bm{f}(\x). \label{eq_poisson}
\end{align}
However, for unboundedly growing vector fields, these numerical integrals diverge.

In \citet{glotzl_2023_helmholtz}, we presented a novel method that generalized the Helmholtz decomposition to $\mathbb{R}^n$ by replacing the usual three-dimensional gradient and vector potentials by an $n \times n$ potential matrix,
allowing to derive gradient and rotation fields:
\begin{defin}
\label{def_potentialmatrix}

Given a vector field $\bm{f} \in C^1(\mathbb{R}^n, \mathbb{R}^n)$, a matrix-valued map $F_{ij} \in C^2(\mathbb{R}^n, \mathbb{R}^{n^2})$ is called \emph{potential matrix} of $\bm{f}$, a scalar function $G \in C^2(\mathbb{R}^{n}, \mathbb{R})$ is called \emph{gradient potential} and a matrix-valued map $R_{ij} \in C^2(\mathbb{R}^n, \mathbb{R}^{n^2})$ is called \emph{rotation potential} of $\bm{f}$, if two conditions are met:
First, $G$ is the sum of the diagonal elements of $F_{ij}$, and $R_{ij}$ is two times the antisymmetric part of $F_{ij}$:
\begin{align}
G(\x) &= \sum_k F_{kk}(\x), \label{eq_def_G} \\
R_{ik}(\x) &= F_{ik}(\x) - F_{ki}(\x). \label{eq_def_R}
\end{align}
Second, the gradient field $\bm{g}$ defined as the gradient of $G$, and the rotation field $\bm{r}$ defined as divergence performed with respect to the rows of the antisymmetric matrix-valued map $R_{ij}$ (here defined as \emph{rotation operator} $\ROT$) yield a Helmholtz decomposition of $\bm{f}$:
\begin{align}
\bm{g}(\x) &= \grad G(\x) = \left[ \partial_{x_i} G(\x);\,_{1 \leq i \leq n}\right] = \left[ \partial_{x_i} \sum\nolimits_k F_{kk}(\x);\,_{1 \leq i \leq n}\right], \label{eq_def_g} \\
\bm{r}(\x) &= \ROT R(\x) \coloneqq \left[ \sum\nolimits_k \partial_{x_k} R_{ik}(\x);\,_{1 \leq i \leq n} \right] = \left[ \sum\nolimits_k \partial_{x_k} \Big( F_{ik}(\x) - F_{ki}(\x) \Big) ;\,_{1 \leq i \leq n} \right], \label{eq_def_r} \\
\bm{f}(\x) &= \bm{g}(\x) + \bm{r}(\x).
\end{align}
\end{defin}
The approach presented in \citet[Sec.~6]{glotzl_2023_helmholtz} to derive such a potential matrix for unboundedly growing fields analytically was quite complicated, requiring distinction of several cases.
With the option to solve Poisson's equation directly following Corollary~\ref{theorem_P}, the approach can be substantially simplified as follows:

\begin{theorem}
\label{theorem_hd}

Let the vector field $\bm f \in C^1(\mathbb{R}^n, \mathbb{R}^n)$ have vector components for which the Poisson equation can be solved such that there exists a vector field $\bm \varphi (\x)$ with $\Delta \bm \varphi (\x) = \bm f(\x)$.
Then, the gradient of $\bm \varphi(\x)$ defined as $F_{ij}(\x) = \partial_{x_j} \varphi_i(\x)$ is a potential matrix of $\bm f$.

\begin{proof}
We know from Lemma~3.4 in \citet{glotzl_2023_helmholtz} (see Appendix~\ref{lemma_hd}) that for every matrix-valued map $F_{ij}$, the gradient field $\bm{g}$ as defined by Eq.~\eqref{eq_def_g} is gradient of some potential, and the rotation field $\bm{r}$ as defined by Eq.~\eqref{eq_def_r} is divergence-free. It remains to be shown that $\bm f(\x) = \bm g(\x) + \bm r(\x)$:
\begin{align}
\bm g(\x) &= \grad G(\x) = \left[ \partial_{x_i} \sum\nolimits_k F_{kk}(\x);\,_{1 \leq i \leq n}\right] = \left[ \partial_{x_i} \sum\nolimits_k \partial_{x_k} \varphi_k(\x);\,_{1 \leq i \leq n}\right], \\
\bm r(\x) &= \left[ \sum\nolimits_k \partial_{x_k} \Big( F_{ik}(\x) - F_{ki}(\x) \Big) ;\,_{1 \leq i \leq n} \right] = \left[ \sum\nolimits_k \partial_{x_k} \Big( \partial_{x_k} \varphi_i(\x) - \partial_{x_i} \varphi_k(\x) \Big) ;\,_{1 \leq i \leq n} \right], \\
\bm f(\x) &= \bm g(\x) + \bm r(\x) = \left[ \partial_{x_i} \sum\nolimits_k \partial_{x_k} \varphi_k(\x)   +    \sum\nolimits_k \partial_{x_k} \Big( \partial_{x_k} \varphi_i(\x) - \partial_{x_i} \varphi_k(\x) \Big)    ;\,_{1 \leq i \leq n}\right] \\
&= \left[ \sum\nolimits_k \partial_{x_k} \partial_{x_k} \varphi_i(\x) ;\,_{1 \leq i \leq n}\right] = \left[ \Delta \varphi_i(\x) ;\,_{1 \leq i \leq n}\right] = \Delta \bm \varphi(\x) = \bm f(\x),
\end{align}
because of the symmetry of second derivatives (Schwarz’s theorem) and interchangeability of sums and derivatives.
\end{proof}
\end{theorem}

Compared to \citet{glotzl_2023_helmholtz}, the approach presented above has three advantages:
First, it is simpler as only $n$ integration problem have to be solved for each vector coordinate, compared to $\tbinom{n}{2} + 1$ for an antisymmetric matrix and a scalar.
Second, it drops the restriction in \citet[Theorem~6.1]{glotzl_2023_helmholtz} stated as ``If $m = k$, then $W(\x)$ must depend only on $x_k$ and $x_m$'' which had no convincing justification anyway.
Third, because we solve the Poisson equation first before differentiating, all integration constants are correctly determined automatically.

The conventional way first differentiates and decomposes at the level of the sources, than integrates to get the potentials, and then differentiates to get the gradient and rotation field.
In \citet{glotzl_2023_helmholtz}, we first jointly differentiated and integrated to get the potential matrix, then decomposed and then differentiated to get the gradient and rotation field. Now, we suggest to first integrate, differentiate once, decompose and differentiate a second time.

\section{Discussion and Conclusions}
\label{sec_conclusion}

This paper introduces the integration-annihilator method to find analytic particular solutions for certain partial differential equations.
In any case, it has to be kept in mind that our methods provide special solutions to inhomogeneous problems.
To get the full solution and to adapt to constraints and boundary conditions, homogeneous solutions have to be added.

The derivations covered various special cases in terms of the differential operators involved and the structure of the inhomogeneous term.
Using examples such as Poisson's equation, the inhomogeneous polyharmonic equation, Helmholtz equation and the wave equation, we demonstrated the usefulness of our approach.
For the Helmholtz decomposition, we showed a simplified approach to fields not decaying at infinity that has a broader applicability than before.
A broad range of inhomogeneities is covered, typically consisting of polynomials, exponential functions, trigonometric functions and their products and sums.
Because polynomials are supported, at least local Taylor approximations can be derived.
A \emph{Mathematica} worksheet implementing the solution strategies can be found in \citet{richters_2025_annihilatorzenodo} and at \url{https://oliver-richters.de/annihilator}.

For every partial differential equation with a linear differential operator, an alternative solution method uses fundamental solutions.
The Malgrange--Ehrenpreis theorem \citep{malgrange_existence_1956, ehrenpreis_solution_1954, ehrenpreis_solution_1955} states that for every linear differential operator with constant coefficients $\D$, a fundamental solution $K$ can be constructed. This distribution is also called Green's function and satisfies $\D K(\x) = \delta(\x)$ where $\delta$ is Dirac's delta function.
For sufficiently smooth functions $q$, a particular solution can be constructed as the convolution $Q(\x) = \int_{\mathbb{R}^n} q(\bm{\xi}) K(\x - \bm{\xi}) d\bm{\xi}^n$ using Green's function as an integral kernel.
Often, this integral needs to be evaluated numerically.
For the kernel integral to be well-defined and finite, the behaviour on $q$ at infinity has to be restricted.
For example, for the case of Poisson's equation, the inhomogeneity must decay faster than $|x|^{-2}$.
With our approach, the restriction on infinite behaviour is not necessary and the particular solution can be derived analytically.

Thanks to their versatility and the possibility of obtaining closed-form solutions, the theorems presented in this paper may prove helpful for problems of (vector) analysis, theoretical physics, and complex systems theory.

\renewcommand{\bibfont}{\normalfont\small}
\addcontentsline{toc}{section}{References} 

\printbibliography

\begin{appendices}

\section{Lemma 3.4 from Glötzl \& Richters (2023)}

\label{lemma_hd}

For every matrix-valued map $F_{ij}$, the gradient field $\bm{g}$ as defined by Eq.~\eqref{eq_def_g} is gradient of some potential, and the rotation field $\bm{r}$ as defined by Eq.~\eqref{eq_def_r} is divergence-free.
If $\bm{f} = \bm{g} + \bm{r}$, the matrix $F_{ij}$ is a `potential matrix of $\bm{f}$'.
\begin{proof}
The gradient field $\bm{g}(\bm{x}) = \grad G(\bm{x})$ is gradient of the potential $G$ as required. The divergence of the rotation field $\bm{r}$ is zero, as
\begin{align}
\begin{split} \Div \bm{r}(\bm{x}) &= \sum\nolimits_i \partial_{x_i} r_i(\bm{x}) = \sum\nolimits_i \partial_{x_i} \sum\nolimits_k \partial_{x_k} R_{ik}(\bm{x}) \\
&= \sum_{i,k} \partial_{x_i} \partial_{x_k} (F_{ik}(\bm{x}) - F_{ki}(\bm{x})) = 0, \end{split} \label{eq_divr}
\end{align}
because the partial derivatives can be exchanged. If $\bm{f} = \bm{g} + \bm{r}$, then $\bm{g}$ and $\bm{r}$ are a `Helmholtz decomposition of $\bm{f}$' as per Definition~\ref{def_helmholtz}, which completes the proof.
\end{proof}

\clearpage

\renewcommand{\thesubsection}{Example \arabic{subsection}}
\section{Examples}

Each example can be derived by following the strategy in Figure~\ref{fig_flowchart}. A \emph{Mathematica} worksheet can be found at \url{https://oliver-richters.de/annihilator} and in \citet{richters_2025_annihilatorzenodo}.

\subsection{Wave equation}
\label{app_wave}

Consider the wave equation that uses $t$ for the temporal dimension of spacetime, and $\Box$ is D'Alembert's operator.
\begin{align}
\Box Q = \left( \partial_t^2/c^2 - \sum_{i=1}^3 \partial_{x_i}^2 \right) Q = q,
\end{align}
with $q(\x) = t \sin(t) x_1^2 x_2$.
Select $\A = \partial_t^2 / c^2$ and $\B = - \sum_{i=1}^3 \partial_{x_i}^2$ and note that for $\lambda = 1$, $\B^{\lambda+1} q(\x) = 0$.
To find a solution for $\A^{\lambda+1} W(\x) = q(\x)$, one needs to find a solution of $\partial_t^4 W / c^4 = q$, a higher-order linear ordinary differential equation in $t$.
A particular solution of $\partial_t^4 f(t) = t \sin(t)$ is known to be $t \sin(t) + 4 \cos(t)$, so selecting $W = (t \sin(t) + 4 \cos(t)) c^4 x_1^2 x_2 $ satisfies all conditions of Corollary~\ref{theorem_WPH} with $k=1$.
This implies that the following special solution satisfies $\Box Q = t \sin(t) x_1^2 x_2$:
\begin{align}
Q &= c^2 x_1^2 x_2 \big(-2 \cos(t) - t \sin(t) \big) + 2 c^4 x_2 \big(4 \cos(t) + t \sin(t) \big). \label{eq_wave_Q}
\end{align}
Calculating $\partial_t^2 / c^2 Q = t \sin(t) x_1^2 x_2 + 2 c^2 x_2 (-2 \cos(t) - t \sin(t))$  and $\Delta Q = 2 c^2 x_2 (-2 \cos(t) - t \sin(t))$ shows that $\Box Q(\x) = t \sin(t) x_1^2 x_2 = q(\x)$ as desired.

Note that the automatic approach in the Mathematica worksheet \citep{richters_2025_annihilatorzenodo} follows Figure~\ref{fig_flowchart} and used $W = c^4 x_1^2 x_2 (t^2+t \sin (t)+4 \cos (t)-4)$ as it integrates $\cos (t^\prime)$ from $0$ to $t$, yielding $\sin (t) - 1$ and therefore additional terms that can be dropped by choosing a different integration constant.
The result $Q = c^2 x_1^2 x_2 (-2 \cos (t) - t \sin (t) + 2) + 2 c^4 x_2 \left(t \sin (t)+4 \cos (t)+t^2-4\right)$ differs from Eq.~\eqref{eq_wave_Q} by a harmonic function.

\subsection{Poisson's equation with multivariate monomial}
\label{app_P_monomial}

Consider Poisson's equation $\Delta^1 Q(\x) = q(\x)$ with a function $q \in C^\infty(\mathbb{R}^2, \mathbb{R})$ with a multivariate monomial given by:
\begin{align}
q(\x) = x_1^2 x_2^{10}.
\end{align}
With $\D = \Delta$ and $k = 1$, Figure~\ref{fig_zusammenhang-groessen} suggests to use $x_2$ as the $W$-coordinate, $\lambda = \lceil (|\beta_{i \neq 2}| -1)/2 \rceil = 1$ with $|\beta_{i \neq 2}| = 2$ the total degree of $x_{i \neq 2}$, $W(\x) = \mathcal{A}_{x_2}^{2\lambda} q(\x) = \frac{10!}{14!} x_1^2 x_2^{14}$ and $\uphi = 1$. With these choices, the conditions of Theorem~\ref{theorem_P} are satisfied, as $\Delta_{\setminus 2}^{\lambda + 1} W(\x) = \partial_{x_1}^4 W(\x) = 0$. The solution of Poisson's equation is given by:
\begin{align}
Q(\x) &= \frac{91 x_1^2 x_2^{12} - x_2^{14}}{12012},
\end{align}
One could also choose $x_1$ as the $W$-coordinate, $\lambda = 5$, $W(\x) = \frac{2!}{14!} x_1^{14} x_2^{10}$ and $\uphi = 1$, yielding 
\begin{align}
Q(\x) &= - \frac{x_1^{14} - 91 x_1^{12} x_2^2 + 1001 x_1^{10} x_2^4 - 3003 x_1^8 x_2^6 + 3003 x_1^6 x_2^8 - 1001 x_1^4 x_2^{10}}{12012},
\end{align}
but the result is less tractable. This shows that the solution is not unique and influenced by the choice of the $W$-coordinate.

\subsection{Polyharmonic equation}
\label{app_PH_sinus}

Consider the polyharmonic equation $\Delta^k Q(\x) = q(\x)$ with a function $q \in C^\infty(\mathbb{R}^2, \mathbb{R})$ given by:
\begin{align}
q(\x) = x_1^2 \sin(x_2).
\end{align}
With $\D = \Delta$, Figure~\ref{fig_zusammenhang-groessen} suggests to use $x_2$ as the $W$-coordinate, $\lambda = \lceil (|\beta_{i \neq 2}| -1)/2 \rceil = 1$ with $|\beta_{i \neq 2}| = 2$ the total degree of $x_{i \neq 2}$, $W(\x) = \mathcal{A}_{x_2}^{2\lambda + 2k} q(\x) = (-1)^{\lambda + k} x_1^2 \sin(x_2)$.
With these choices, the conditions of Corollary~\ref{theorem_PH} are satisfied, as $\Delta_{\setminus 2}^{\lambda + 1} W(\x) = \partial_{x_1}^4 W(\x) = 0$. The solution of the polyharmonic equation is given by:
\begin{align}
Q(\x) &= \sum_{p=0}^{\lambda} (-1)^p \tbinom{k+p-1}{p} \partial_{x_1}^{2\lambda - 2p} \partial_{x_2}^{2p} W(\x) \\
&= \partial_{x_2}^{2} W(\x) - \tbinom{k}{1} \partial_{x_1}^2 W(\x) \\
&= (-1)^k x_1^2 \sin(x_2) + k (-1)^k \cdot 2 \sin(x_2) \\
&= (-1)^k \sin(x_2) \left(x_1^2 + 2 k \right).
\end{align}
One can show by simple calculation that applying the Laplacian $d$-times results in:
\begin{align}
\Delta^d Q(\x) &= (-1)^d (-1)^k \sin(x_2) \left(x_1^2 + 2 k - 2 d \right), \\
\Rightarrow \Delta^k Q(\x) &= x_1^2 \sin(x_2).
\end{align}

\label{sec:last_example}

\end{appendices}

\end{document}